\documentclass[12pt,twoside]{article}
\usepackage{fullpage}

\usepackage{multicol}

\usepackage{amsmath}
\usepackage{amssymb}

\numberwithin{equation}{section}
\newcommand{\spec}{\mathrm{Spec}}
\renewcommand{\star}{\mathrm{star}}
\renewcommand{\P}{\mathbb{P}}
\newcommand{\E}{\mathbb{E}}

\newtheorem{thm}{Theorem}[section]
\newtheorem{defi}[thm]{Definition}

\newtheorem{lem}[thm]{Lemma}

\newtheorem{rem}[thm]{Remark}

\newcommand{\spr}{\mathrm{spr}}

\begin{document}


\title{Ramanujan Graphs and Interlacing Families}
\date{August 20, 2024}

\author{Nikhil Srivastava}

\maketitle
\begin{abstract}
This survey accompanies a lecture on the paper ``Interlacing Families I: Bipartite Ramanujan Graphs of All Degrees'' \cite{marcus2015interlacing} by A. Marcus, D. Spielman, and N. Srivastava which was awarded a Frontiers of Science Award at the 2024 ICBS in July, 2024. Its purpose is to explain the developments surrounding this work over the past ten or so years, with an emphasis on connections to other areas of mathematics. Earlier surveys about the interlacing families method by the same authors \cite{marcus2014ramanujan,marcus2016solution} focused on applications in functional analysis, whereas the focus here is on applications in spectral graph theory.
\end{abstract}

\tableofcontents

\section{Introduction} 
Ramanujan graphs are sparse graphs with spectral gap as
large as possible. These graphs were introduced by Lubotzky-Phillips-Sarnak \cite{lubotzky1988ramanujan} and Margulis \cite{margulis1988explicit} and have remarkable
connections to many areas of mathematics and computer science, including
spectral graph theory, random matrices, number theory, pseudorandomness, metric
geometry, free probability, and the geometry of polynomials. In this article I
will explain some relatively recent constructions of Ramanujan graphs based on the ``method of interlacing families of polynomials'' and discuss the new connections they reveal. 

The study of Ramanujan graphs can be motivated in at least two different ways.

\subsection{Spectral Expansion}
One context in which Ramanujan graphs is spectral graph theory. Let us review the basic setup. Let $G$ be a $d-$regular undirected graph on $n$ vertices with adjacency matrix $A$. Label the eigenvalues of $A$ in descending order: $$\lambda_1(A)\ge \lambda_2(A)\ge\ldots\ge \lambda_n(A).$$ Recall\footnote{If the reader has no memory of such statements, then it is a nice exercise to prove them.} that:
\begin{itemize}
    \item We always have $\lambda_1(A)=d$, since the row sums of $A$ are equal to the degree of $G$.
    \item $\lambda_1(A)>\lambda_2(A)$ if and only if $G$ is a connected graph.
    \item $\lambda_n(A)=-d$ if and only if $G$ is a bipartite graph, and in this case the spectrum of $A$ is symmetric about zero. If $\lambda_n(A)=-d$ and $G$ is connected then this eigenvalue must accordingly be simple.
\end{itemize}
We will accordingly refer to $\lambda_i,i\neq 1$ as the {\em nontrivial} eigenvalues of $A$ when $G$ is nonbipartite, and to $\lambda_i,i\neq 1,n$ as the nontrivial eigenvalues of $A$ when $G$ is bipartite. We can now state the definition of a Ramanujan graph.
\begin{defi}\label{def:raman}
    A $d-$regular graph $G$ is {\em Ramanujan} if it is connected and $$\spec(A)\subset [-2\sqrt{d-1},2\sqrt{d-1}]\cup\{-d,d\},$$
    i.e.,  all of its nontrivial eigenvalues lie in $[-2\sqrt{d-1},2\sqrt{d-1}]$.
\end{defi}
It is easy to find graphs of small size which are Ramanujan; both $K_d$ and $K_{d,d}$ are examples. The challenge is to construct infinite sequences which are Ramanujan. The significance of the number $2\sqrt{d-1}$ comes from the following lower bound on the size of nontrivial eigenvalues.
\begin{thm}[Alon-Boppana \cite{alon1987monotone}]\label{thm:ab} Suppose $G$ is a $d-$regular graph on $n$ vertices. Then there is a nontrivial eigenvalue of its adjacency matrix $A$ satisfying:
\begin{equation}\label{eqn:alonboppana} |\lambda_i(A)|\ge 2\sqrt{d-1}-o_n(1).\end{equation}
\end{thm}
Thus, up to the $o_n(1)$ term in the Alon-Boppana bound, Ramanujan graphs have the largest spectral gap among all $d-$regular graphs.  

The number $2\sqrt{d-1}$ is equal to the spectral radius of the adjacency
operator $A_T$ of the infinite $d-$ary tree $T_d$, in which
all vertices have degree $d$.
The Alon-Boppana bound can be proven in 
at least two ways: by showing that (generalized) eigenfunctions of $A_{T_d}$ can be truncated to produce test functions 
on $A_G$ with appropriately large or small quadratic form, or by relating the moments of the eigenvalues of $A_G$ to the
moments of a certain spectral measure of $A_T$ by observing that $G$ has at least as many closed walks of a certain
length as $T$ does.\\

\noindent {\em Covers.}  The significance of $T_d$ is that it is the {\em universal cover} of
every $d-$regular graph.  Recall that a graph $H=(V_H,E_H)$ is a cover of a graph $G=(V_G,E_G)$ if there
is a graph homomorphism\footnote{ i.e., a map
satisfying $\pi(x)\pi(y)\in E_G$ whenever $xy\in H$.} $\pi:V_H\rightarrow V_G$ such that for every $x\in
V_G$, the preimage $\pi^{-1}(\star_G(x))$ is a disjoint union of stars\footnote{$\star_G(x)$ is the 
induced subgraph of $G$ on $x$ and its neighbors.} in $H$ which are isomorphic to $\star_G(x)$. The number of 
stars in the preimage is called the {\em degree} of the cover. For a given $d-$regular graph $G$, the infinite tree $T_d$ is called its universal cover because it covers every finite graph which covers $G$.  This property is
used to compare eigenfunctions and walk counts between $G$ and $T_d$ in the proofs of Theorem \ref{thm:ab}.

\begin{rem} There is a generalization of the Alon-Boppana bound due to Greenberg \cite{greenberg1995spectrum} , who
showed that for any (not necessarily regular) sequence of graphs $G_n$ on $n$ vertices with a
common universal covering tree $T$, there is a nontrivial eigenvalue of $G_n$ satisfying:
$$|\lambda_i(G_n)|\ge \spr(A_T)-o_n(1),$$
where $A_T$ is the adjacency operator of $T$ and $\spr(\cdot)$ denotes the spectral radius. We will only discuss regular
graphs in this survey for simplicity, but the results generalize to the irregular case.
\end{rem}

The most important result about Ramanujan graphs is that they exist. This was shown independently by Margulis \cite{margulis1988explicit} and Lubotzky-Phillips-Sarnak \cite{lubotzky1988ramanujan} using deep results from number theory, which are responsible for the Ramanujan nomenclature.
\begin{thm}[\cite{lubotzky1988ramanujan,margulis1988explicit}]
    Let $d=p^k+1$. Then there are infinite sequences of both nonbipartite and bipartite $d-$regular Ramanujan graphs.
\end{thm}
This construction left open the question of what happens for other degrees.\\

\noindent {\bf Question 1 (Lubotzky \cite{lubotzky1994discrete}).} Do infinite sequences of Ramanujan graphs exist for all degrees $d\ge 3$?\\

Another question is whether such graphs exist for all values of $n$, subject to the obvious constraints imposed by regularity. In this survey we will describe partial answers to these questions.

\subsection{Random Graphs}\label{sec:introrand}
There is a second, probabilistic, context in which Ramanujan graphs appear naturally: random graphs. We will consider two models of random graphs:
\begin{itemize}
    \item[(R1)] Given $d\ge 3$ and $n$ even, choose $G$ uniformly at random from the set of all $d-$regular graphs on $n$ vertices. We will also consider the related model $(R1')$ in which $G$ is a union of $d$ random perfect matchings  on $[n]$ (assuming $n$ is even). These models are known to be contiguous, in the sense that any high probability event in one model also holds with high probability in the other \cite{wormald1999models}. Note that the second model produces graphs which may have multiple edges.
    \item[(R2)] Given a fixed $d-$regular base graph $H$ and $n\ge 1$, choose a uniformly random $n-$cover $G$ of $H$. In this case $H$ is a $d-$regular graph on $n|V_H|$ vertices.
\end{itemize}
Alon conjectured that a random graph from model (R1) should typically be ``almost''
Ramanujan, i.e., its adjacency matrix should satisfy
\begin{equation}\label{eqn:alonconj} |\lambda_i|\le 2\sqrt{d-1}+\epsilon_n\qquad
\forall i\neq 1,\end{equation} with probability $1-o_n(1)$ where
$\epsilon_n=o_n(1)$. Alon's conjecture was proven by Friedman \cite{friedman2008proof} with
$\epsilon_n=O_d((\log\log n/\log(n))^2)$ (see also a simpler proof by Bordenave
\cite{bordenave2015new}, and another proof by Puder \cite{puder2015expansion} which obtains $\epsilon_n=1$). Two very recent works
significantly improve this. Huang, McKenzie, and Yau \cite{huang2024optimal} (building on \cite{bauerschmidt2020edge,huang2022edge}) obtain $\epsilon_n=n^{-2/3+o_n(1)}$, which
is conjectured to be sharp as it is believed that the fluctuations of $\lambda_2$ should (after appropriate normalization) behave like the fluctuations of a random GUE matrix (i.e., obey Tracy-Widom statistics).
Chen et al. \cite{chen2024new} obtain $\epsilon=n^{-1/8+o_n(1)}$ by a simpler method in the more
general context of strong convergence in free probability theory (see also the influential work of Bordenave-Collins \cite{bordenave2019eigenvalues}, who introduced this connection). 

Friedman generalized \cite{friedman2003relative} Alon's conjecture to the model (R2), asking whether
\eqref{eqn:alonconj} continues to hold in this setting, where one restricts attention to the ``new'' eigenvalues of $A_G$ (see Section \ref{sec:randcovers} for a definition). Friedman and Kohler \cite{friedman2019relativized}
proved the generalized conjecture with $\epsilon_n=o_n(1)$, and Puder \cite{puder2015expansion} proved it up to a multiplicative
context in the general setting of possibly irregular base graphs $H$, where
$2\sqrt{d-1}$ is replaced by $\spr(A_T)$ for $T$ the universal cover of $H$. Different proofs of the generalized Alon conjecture were given  in the context of strong convergence by Bordenave-Collins \cite{bordenave2019eigenvalues} and Chen et al. \cite{chen2024new}.

Given the above background, a natural question is the following:\\

\noindent {\bf Question 2.} Can we obtain $\epsilon_n=0$ in the above bounds if we
settle for a lower probability bound?

Such a result has been conjectured by Miller, Novikov, and Sabelli \cite{miller2008distribution}, based on numerical experiments. It is known that this cannot hold with probability $1-o_n(1)$ due
to the possibility of certain constant sized subgraphs appearing in $G$ which can cause $\lambda_2(G)>2\sqrt{d-1}$. This presents an obstacle to employing the standard methods of random matrix theory, which typically yield bounds which hold with high probability, a conclusion which is too strong in this context.
\section{Interlacing Families}
In this section we describe a method which is able to establish {\em low probability}
control on the eigenvalues of certain random matrices; in fact, it simply says
that a certain probability is nonzero, without providing any estimate.

Given a random $n\times n$ Hermitian matrix $A$, we can consider its {\em expected characteristic polynomial}:
$$p_A(z):=\E \det(zI-A).$$
Expected characteristic polynomials of random unitary matrices \cite{diaconis1994eigenvalues, keating2000random} and Gaussian matrices \cite{diaconis2004random} were studied in the context of number theory, where they are believed to be related to the zeros of the Riemann zeta function. Expected characteristic polynomials of certain Hermitian matrices were studied by Godsil and Gutman \cite{godsil1978matching}, in the context of algebraic graph theory.

The works \cite{marcus2015interlacing, marcus2018interlacing} discovered a new method relating the zeros of the polynomial $p_A(z)$ to the eigenvalues of the random matrix $A$ under certain conditions. This relationship is nontrivial since the maps relating a matrix, its characteristic polynomial, and its eigenvalues are highly nonlinear. It is useful because it is often much easier to compute $p_A$ then to compute statistics of the eigenvalues of $A$. The following theorem follows from results in \cite{marcus2015interlacing, marcus2018interlacing, hall2018ramanujan}.
\begin{thm}[Interlacing Families] \label{thm:if} Let $A$ be the adjacency matrix of a random graph chosen from model $(R1')$ or $(R2)$ as described in Section \ref{sec:introrand}. Then $p_A(z)$ is real-rooted and for every $i\le n$:
$$ \P[ \lambda_i(A)\le \lambda_i(p_A)]>0,$$
where we use $\lambda_i$ to denote the $i$th largest root of a real-rooted polynomial.
\end{thm}
We refer the reader to \cite{marcus2015interlacing,marcus2018interlacing, hall2018ramanujan} for a proof of this theorem.  The key mechanism behind the theorem is encapsulated in the following lemma about a pair of real-rooted univariate polynomials, which the interested reader may prove as an exercise.
\begin{lem}\label{lem:fell} Suppose $p_0(z)$ and $p_1(z)$ are monic polynomials of degree $n$ such that all of the convex combinations
$$ p_t:= (1-t)p_0+tp_1,\quad t\in [0,1]$$
are real-rooted. Then, for every $i$ and $t\in [0,1]$:
$$ \min\{\lambda_i(p_0),\lambda_i(p_1)\}\le \lambda_i(p_t)\le \max\{\lambda_i(p_0),\lambda_i(p_1)\}.$$
\end{lem}
In particular, Lemma \ref{lem:fell} allows one to relate the roots of an average of two polynomials to the roots of the individual polynomials. Theorem \ref{thm:if} is proven by writing $p_A$, an expectation over exponentially many polynomials, as a repeated conditional expectation over pairs of polynomials, and inductively applying Lemma \ref{lem:fell}. The required real-rootedness statements are established using the theory of multivariate real stable/hyperbolic polynomials, especially the work of Borcea-Branden \cite{borcea2009lee,borcea2009lee2}. The term {\em interlacing family} refers to this inductive scheme. 

Theorem \ref{thm:if} has been generalized to matrix models other than random graphs --- see e.g. \cite{if2, ravichandran2020mixed,anari2014kadison,bownik2023akemann} for other $A$ to which the theorem applies. The common feature of all of these models is that it is possible to generate the random matrix $A$ by a sequence of independent low rank (rank $1$ or $2$ perturbations of a fixed matrix). It is possible Theorem \ref{thm:if} generalizes to other models which have not yet been studied. It is also easy to see that it does not hold for arbitrary random matrices $A$, even in $n=2$.

\newcommand{\bk}[2]{|#1\rangle \langle #2|}
\newcommand{\C}{\mathbb{C}}
\section{Random Covers of Fixed Graphs}\label{sec:randcovers}
In this section we explain how Theorem \ref{thm:if} applied to the random model $(R2)$ can be used to partially answer Questions $1$ and $2$. Fix a base graph $H$ and let $G$ be a random $n-$cover of $H$. Concretely, the adjacency matrix $A_G$ of $G$ is given by:
\begin{equation}\label{eqn:liftadj}
A_G = \sum_{ij\in E_H} \bk{i}{j}\otimes P_{ij}+\bk{j}{i}\otimes P_{ij}^*,
\end{equation}
where $\{P_{ij}\}_{ij\in E_H}$ are independent random $n\times n$ permutation matrices indexed by the edges of $H$. Notice that if $v$ is an eigenvector of $A_H$ satisfying $A_H v = \lambda v$, then 
$$ A_G (v\otimes \mathbf{1}) = (A_H v)\otimes \mathbf{1}=\lambda (v\otimes \mathbf{1}),$$
since $\mathbf{1}\in\C^n$ (the all ones vector) satisfies $P_{ij}\mathbf{1}=P_{ij}^*\mathbf{1}=\mathbf{1}$ for all $ij\in E_H$. Thus, the spectrum of $A_G$ contains the spectrum of $A_H$. Following Friedman \cite{friedman1993some}, we will refer to $\spec(A_G)\setminus \spec(A_H)$ (counting eigenvalues with multiplicity) as the ``new'' eigenvalues of $A_G$ relative to $H$.

\subsection{$2-$covers} Bilu and Linial \cite{bilu2006lifts} studied the case of $2-$covers. They observed that in this case, the new eigenvalues of $A_G$ for a given cover (defined by permutations $P_{ij}$) are equal to the eigenvalues of the corresponding {\em signed adjacency matrix} of $A_H$:
$$ A_{H,s}(i,j):= \det(P_{ij}),\quad ij\in E_H.$$
They showed that every $H$ has a signing such that $\|A_{H,s}\|\le O(\sqrt{d\log^3 d})$ and conjectured that this bound should be improvable to $2\sqrt{d-1}$. This conjecture would yield a way to construct infinite sequences of (bipartite and nonbipartite) Ramanujan graphs by starting with $K_d$ or $K_{d,d}$ as a base graph and repeatedly finding $2-$covers whose new eigenvalues satisfy the Ramanujan bound.

Bilu and Linial's conjecture is still open, but \cite{marcus2015interlacing} showed that it is true for bipartite graphs.
\begin{thm}\label{thm:if1}
    If $H$ is a $d-$regular bipartite graph, then there exists a signing $A_{H,s}$ such that $\|A_{H,s}\|\le 2\sqrt{d-1}$.
\end{thm}
The proof of this theorem has three ingredients. The first is Theorem \ref{thm:if}. The second is a combinatorial expression for the expected characteristic polynomial of a random signing of a graph $G$ in terms of its {\em matching polynomial} $M_G$, discovered by Godsil and Gutman in 1978 \cite{godsil1978matching}.
\begin{thm}[Godsil-Gutman]\label{thm:gg} Let $H$ be an undirected graph and let $A_{H,s}$ be a uniformly random signing of its adjacency matrix. Then,
$$ \E \det(zI-A_{H,s})=\sum_{k=0^\lfloor n/2\rfloor} x^{n-2k}(-1)^k m_k(H)=:M_H(z),$$
where $m_k(H)$ is the number of matchings with $k$ edges in $H$.
\end{thm}

The third ingredient comes from statistical mechanics, in which the matching polynomial $M_H$ was studied a decade earlier by Heilmann and Lieb \cite{heilmann1972theory} in the context of thermodynamic phase transitions in the monomer dimer model, where it arises as a partition function.
\begin{thm}[Heilmann-Lieb]\label{thm:hl} Let $H$ be an undirected graph with maximum degree $d$. Then, the matching polynomial $M_H(z)$ is real-rooted and all of its roots are contained in $[-2\sqrt{d-1},2\sqrt{d-1}]$.
\end{thm}
The appearance of the number $2\sqrt{d-1}$ in Theorem \ref{thm:hl} is not an accident: it is proven by an inductive argument which deletes the vertices of $H$ one at a time, whose recursion tree is a subtree of $T_d$. As 
 observed in \cite{marcus2015interlacing}, this argument can also be generalized to universal covers of irregular graphs.

Given these three ingredients, it follows immediately that $\lambda_1(A_{H,s})\le 2\sqrt{d-1}$ with nonzero probability. Since bipartite graphs have symmetric spectrum, this yields that $\|A_{H,s}\|\le 2\sqrt{d-1}$, which is the only way in which bipartiteness is used in the proof.
\subsection{$n-$covers} Hall, Puder, and Sawin \cite{hall2018ramanujan} generalized the results of the previous subsection in two ways. First, they showed that the same thing is true if one considers $n-$covers (as defined in \eqref{eqn:liftadj}) instead of $2-$covers.
\begin{thm}[\cite{hall2018ramanujan}]\label{thm:hps1}
    If $H$ is a $d-$regular bipartite graph and $n\ge 2$, then there exists an $n-$cover $G$ of $H$ such that the new eigenvalues of $A_G$ are bounded in absolute value by $2\sqrt{d-1}$.
\end{thm}
Second, they considered ``signings'' of a base graph $H$ by matrices arising from a group representation. Let $\Gamma$ be a group and let $\rho:\Gamma\rightarrow GL_n(\C)$ be a unitary representation of $\gamma$. Given a labeling $s:E_H\rightarrow \Gamma$, consider the operator
$$A_{H,s} = \sum_{e\in E_H} \bk{i}{j}\otimes \rho(s(e))+\bk{j}{i}\otimes \rho(s(e))^*,$$
which we will call a {\em signing of $H$ by $\rho$}.
In the special case of $\Gamma=S_n$ with $\rho$ its regular representation, the eigenvalues of $A_{H,s}$ are the eigenvalues of the adjacency matrix of the cover defined by $s$. The ''old'' eigenvalues arise from the trivial representation of $S_n$ and the new eigenvalues arise from the remaining representations. 
\begin{thm}[\cite{hall2018ramanujan}]\label{thm:hps2}
    Suppose $\Gamma$ is a group and $\rho:\Gamma\rightarrow GL_n(\C)$ is a unitary representation with two properties:
    \begin{itemize}
        \item [(P1)] All exterior powers $\bigwedge^m \rho$, $0\le m\le n$, are irreducible and non-isomorphic.
        \item [(P2)] $\rho(\Gamma)$ is a complex pseudo-reflection group, i.e., it is generated by matrices $A$ which satisfy
        $rank(A-I)=1$.
    \end{itemize}
    Then, if $H$ is a bipartite $d-$regular graph, there is a signing $A_{H,s}$ of $H$ by $\rho$ satisfying $\|A_{H,s}\|\le 2\sqrt{d-1}$.
\end{thm}
The proof of Theorem \ref{thm:hps2} relies on two ingredients. The first is a generalization of Theorem \ref{thm:if} to the setting of signings by representations, showing that
$$ \P[ \lambda_1(A_{H,s})\le \lambda_1(\E \det(zI-A_{H,s})]>0.$$
Here, property $(P2)$ is used to produce an interlacing family. The second ingredient is a generalization of the Godsil-Gutman theorem, showing that if $\rho:\Gamma\rightarrow GL_n(\C)$ satisfies $(P1)$ then:
\begin{equation}\label{eqn:p2avg} \E \det(zI - A_{H,s}) = \E_{G\sim Cov_{n-1}} M_G(z),\end{equation}
where $G$ is a uniformly random $n-1$-cover of $H$ and $M_G$ is the matching polynomial as in Theorem \ref{thm:gg}. Interestingly, the expected characteristic polynomial ``forgets'' almost everything about $\rho$. The property $(P1)$ is used to invoke the Peter-Weyl theorem in a way which produces the required cancellations. The identity \eqref{eqn:p2avg} implies that $\lambda_1(\E \det(zI-A_{H,s}))\le 2\sqrt{d-1}$ since this is true for each $M_G$ in the average by Heilmann-Lieb. Together with the generalization of Theorem \ref{thm:if}, the proof is concluded.
\subsection{Beyond Trees}
All of our discussion so far has been about comparing the spectrum of a finite graph with that of its universal cover $T$, which is an infinite tree. But the same kind of question makes sense even when the infinite covering graph is not a tree. Mohanty and O'Donnell \cite{mohanty2019x} generalized the results of Hall-Puder-Sawin to the setting of certain infinite graphs called {\em additive product graphs} (see their paper for a definition) and proved the following.
\begin{thm}\label{thm:mo} Suppose $X$ is an infinite additive product graph and $H$ is a finite graph covered by $X$. Then for every $n\ge 2$ there is an $n-$cover $G$ of $H$ such that every new eigenvalue $\lambda\in \spec(A_G)\setminus\spec(A_H)$ satisfies $\lambda\le \spr(A_X)$.
\end{thm}
The class of additive product graphs includes several interesting examples which are not trees, notably certain Cayley graphs of the modular group $PSL_2(\mathbb{Z})$.

The proof of Theorem \ref{thm:mo} requires analyzing a new family of expected characteristic polynomials which they refer to as ``additive characteristic polynomials'' (see also \cite{ravichandran2020mixed}), which includes all of the other expected characteristic polynomials considered in this survey as special cases. A remarkable step in the proof is to prove a root bound for all such polynomials using algebraic combinatorics, specifically Viennot's theory of heaps of pieces \cite{viennot2006heaps}, which may be viewed as generalizing Heilmann and Lieb's tree recursion proof beyond trees. 

Theorem \ref{thm:mo} makes progress towards a question originally asked by Clark \cite{clark2006ramanujan}: given an infinite graph $X$ with a finite quotient $H$ satisfying 
\begin{equation}\label{eqn:xram} \lambda_2(A_H)\le \spr(A_X)\end{equation}
must it have infinitely many quotients satisfying \ref{eqn:xram}? The condition on the existence of a finite quotient is required, as Lubotzky and Nagnibeda \cite{lubotzky1998not} showed that there are graphs where no finite quotient has this property.

\section{Random Regular Graphs}
In this section we consider the random matrix model $(R1')$ (random $d-$regular graphs) defined in Section \ref{sec:introrand}. Let $n$ be even and let $A_M$ be the adjacency matrix of a perfect matching on $n$ vertices.
Then the adjacency matrix of a random graph sampled from the model $(R1')$ is given by:
$$ A_G = \sum_{i\le d} P_i A_M P_i^*,$$
where $P_1,\ldots,P_d$ are independent random $n\times n$ permutations. We then have:
\begin{thm}[\cite{marcus2018interlacing}]\label{eqn:if4result} For $A_G$ as above, $$\P[\lambda_2(A_G)\le 2\sqrt{d-1}]>0.$$\end{thm}
i.e., $A_G$ is ``one-sided Ramanujan'' with nonzero probability\footnote{Actually, \cite{marcus2018interlacing} showed that the same is true for random bipartite $d-$regular graphs, which are consequently ``two-sided'' Ramanujan per Definition \ref{def:raman}, but we stick to the one-sided non-bipartite case in this section as its proof is simpler.}.
By Theorem \ref{thm:if}, this conclusion follows from the bound 
\begin{equation}\label{eqn:if4bound} \lambda_2(\E \det(zI-A_G))\le 2\sqrt{d-1},\end{equation}
which we will focus on for the remainder of this section.

The expected characteristic polynomial in \eqref{eqn:if4bound} is not a matching polynomial, but it turns out to have another special structure. In particular, it can be related in a simple way to 
$$\chi[A_M](z):=\det(zI-A_M)=(z^2-1)^{n/2},$$ 
the characteristic polynomial of a single matching. 
The relation is via the {\em Walsh Convolution} of two polynomials, denoted $\boxplus_n$, which defined as follows. Given any monic polynomial $p(z)$ of degree $n$, let $\hat{p}$ be the unique polynomial of degree $n$ satisfying
\begin{equation}\label{eqn:fourierdef} \hat{p}(d/dz)z^n = p(z). 
\end{equation}
Given a pair of monic polynomials $p$ and $q$ of degree $n$, define
$$ (p \boxplus_n q)(z) := \hat{p}(d/dz)\hat{q}(d/dz) z^n.$$
This convolution was originally discovered by Walsh in \cite{walsh1922location} in the context of classical invariant theory, where he showed that it preserves real-rootedness. Its connection to expected characteristic polynomials of random matrices was discovered in \cite{marcus2022finite}. 

The punch line is that for $p(z):=\chi[M](z)/(z-1) = (z-1)^{n/2-1}(z+1)^{n/2}$ (the characteristic polynomial of a matching with the trivial root removed), we have:
\begin{equation}\label{eqn:fflinearize}
   \E \det(zI-A_G)/(z-d) = p\boxplus_n p \boxplus_n\ldots\boxplus_n p \quad\textrm{$d$ times},
\end{equation}
i.e., once we drop the trivial eigenvalue, the characteristic polynomial of a sum of random matchings is equal to a convolution of characteristic polynomials of single matchings. The proof of this identity appears in \cite{marcus2018interlacing} and boils down to showing that the polynomial $(A_1,\ldots,A_d)\mapsto \E\det(zI-\sum_{i\le d} P_i A_iP_i^*)$ (for Hermitian matrices $A_i$) is invariant under conjugating the $A_i$ by unitary matrices which fix the $\mathbf{1}$ vector. It is also shown in \cite{hall2018ramanujan} using elementary representation theory.

\newcommand{\R}{\mathbb{R}}
Perhaps even more surprisingly, the locations of the roots of a Walsh convolution of polynomials can be related to another convolution appearing in a different area of mathematics: Voiculescu's {\em free convolution} of probability measures \cite{voiculescu1991limit}. We may define (a special case of) it in the present context as follows: given two discrete measures on $\R$:
$$\mu = (1/n)\sum_{i\le n} \delta_{x_i}, \quad \nu = (1/n)\sum_{i\le n}\delta_{y_i},$$
let $A_N, B_N$ be a sequence of $Nn\times Nn$ diagonal matrices with eigenvalues $\{x_i\}_{i\le n},\{y_i\}_{i\le n}$ (respectively), each with multiplicity $N$. Then the empirical spectral distribution of $$A_N+Q_NB_NQ_N^*$$ where $Q_N$ is a Haar unitary converges weakly to a compactly supported measure on $\R$. The free convolution $\mu\boxplus \nu$ is defined to be this measure.

We then have the following relation between the Walsh convolution and Voiculescu's convolution, essentially\footnote{\cite{marcus2022finite} proves a certain inequality of Cauchy transforms from which Theorem \ref{thm:ffc} follows quickly, see Appendix B of \cite{mohanty2019x} for a discussion. } shown in \cite{marcus2022finite}.
\begin{thm}\label{thm:ffc}
    If $p_1,\ldots,p_d$ are real-rooted polynomials of degree $n$ and $\mu_1,\ldots,\mu_d$ are the uniform measures on their roots, then
    $$ \lambda_1(p_1\boxplus_n p_2\boxplus_n\ldots\boxplus_n p_d)\le \max\mathrm{supp}(\mu_1\boxplus\mu_2\boxplus\ldots\boxplus \mu_d).$$
\end{thm}
Applying this theorem\footnote{along with a monotonicity argument to replace the uniform measure on the roots of $p$ by the uniform measure on the roots of $\chi[M](z)$} 
to the right hand side of $\eqref{eqn:fflinearize}$, we have
$$\lambda_1(p\boxplus_n\ldots\boxplus_n p)\le \max\mathrm{supp}(\mu\boxplus\ldots\boxplus\mu),$$
where $\mu=(1/2)(\delta_{-1}+\delta_{+1})$. The free convolution of $d$ Bernoulli measures $\mu\boxplus\ldots\mu$ is equal to the Kesten-Mckay law \cite{mckay1981expected}, the limiting empirical spectral distribution of random $d-$regular graphs, which is supported on $[-2\sqrt{d-1},2\sqrt{d-1}]$. This completes the proof of Theorem \ref{eqn:if4result}. 

There are three significant differences between Theorem \ref{eqn:if4result} (and its bipartite version \cite{marcus2018interlacing}) and Theorem \ref{thm:if1}: (i) It obtains multigraphs instead of simple graphs. (ii) It produces (bipartite) Ramanujan graphs of every size, instead of just an infinite sequence of them. (iii) Cohen \cite{cohen2016ramanujan} has shown that the graphs guaranteed by Theorem \ref{eqn:if4result} can be found in polynomial time via an ingenious alternate interlacing family which can be efficiently computed. In contrast, we do not know how to find Ramanujan $2-$covers or $n-$covers of arbitrary base graphs in polynomial time, and the relevant expected characteristic polynomials are $\#P$-hard to compute exactly (though their roots can be approximated up to a constant factor in subexponential time, see \cite{anari2018approximating}).

The methods of this section seem to be quite different from those in Section 3, relying on ideas from invariant theory and free probability rather than algebraic combinatorics and statistical mechanics. It is worth noting that there is an alternate proof of Theorem \ref{eqn:if4result} following the approach of \cite{hall2018ramanujan}. The idea is to view a graph from the model $(R1')$ as a random $n-$covering of a ``bouquet''  (a single-vertex multigraph with $d$ self-loops). Hall-Puder-Sawin conjectured that their proof should work in this setting, but could only show that it works for graphs without self-loops (importantly, they showed that it works for covers of a $2-$vertex graph with $d$ parallel edges, which corresponds to a random bipartite $d-$regular graph). Their conjecture was proven by Amini \cite{amini2019stable} using ideas from stability theory.

Finally, the study of expected characteristic polynomials of sums and products of unitarily invariant random matrices led to a theory of ``finite free probability'' which is currently an active research topic (see e.g. \cite{marcus2021polynomial,marcus2022finite,arizmendi2018cumulants,arizmendi2023finite,arizmendi2024s,gorin2020crystallization,gribinski2024theory}).
\section{Open Questions}
The common theme in the works discussed in this survey is:
\begin{quote}
    Expected characteristic polynomials allow us to {\em sharply} relate the spectra of finite dimensional random matrices and
    the infinite dimensional operators to which they converge, but with {\em low probability}.
\end{quote}
There is likely more to be said about this connection. We conclude with six questions which we hope will lead to further developments.
\begin{enumerate}
    \item Are there infinite sequences of non-bipartite Ramanujan graphs for every $d\ge 3$? The above theorems only handle the bipartite case.
    \item Is the probability that a random regular graph is Ramanujan bounded away from zero?
    \item Is there a polynomial time algorithm for finding the covers guaranteed by Theorems \ref{thm:if1} and \ref{thm:hps1}? Currently, such an algorithm is only known in the setting of Theorem \ref{eqn:if4result} \cite{cohen2016ramanujan}.
    \item More speculatively, can the methods discussed above be employed in the setting of random hyperbolic surfaces? The recent breakthrough of Hide and Magee \cite{hide2023near} shows that random covers of such surfaces are the appropriate analogue of ``Ramanujan'' in that setting.
    \item Can the results of this survey be extended to obtain ``super-Ramanujan'' bounds of type
    $$ |\lambda_i(A_G)|\le 2\sqrt{d-1}-\epsilon_n,$$
    for some $\epsilon_n$ asymptotically matching the Alon-Boppana bound? The proofs above yield roughly $\epsilon_n=d/n$, but it is conceivable that $1/\log^c(n)$ is achievable. This question was asked by P. Sarnak in the context of computing the ``bass note spectrum'' of various homogeneous spaces.
    \item Following \cite{hall2018ramanujan}, does every $d-$regular bipartite graph have a Ramanujan signing by elements of $PSL_2(\mathbb{F}_q)$?
\end{enumerate}
\bibliographystyle{plain}
\bibliography{icbs}

\begin{thebibliography}{10}

\bibitem{alon1987monotone}
Noga Alon and Ravi~B Boppana.
\newblock The monotone circuit complexity of boolean functions.
\newblock {\em Combinatorica}, 7:1--22, 1987.

\bibitem{amini2019stable}
Nima Amini.
\newblock Stable multivariate generalizations of matching polynomials.
\newblock {\em arXiv preprint arXiv:1905.02264}, 2019.

\bibitem{anari2014kadison}
Nima Anari and Shayan~Oveis Gharan.
\newblock The kadison-singer problem for strongly rayleigh measures and
  applications to asymmetric tsp.
\newblock {\em arXiv preprint arXiv:1412.1143}, 2014.

\bibitem{anari2018approximating}
Nima Anari, Shayan~Oveis Gharan, Amin Saberi, and Nikhil Srivastava.
\newblock Approximating the largest root and applications to interlacing
  families.
\newblock In {\em Proceedings of the Twenty-Ninth Annual ACM-SIAM Symposium on
  Discrete Algorithms}, pages 1015--1028. SIAM, 2018.

\bibitem{arizmendi2024s}
Octavio Arizmendi, Katsuori Fujie, Daniel Perales, and Yuki Ueda.
\newblock $ s $-transform in finite free probability.
\newblock {\em arXiv preprint arXiv:2408.09337}, 2024.

\bibitem{arizmendi2023finite}
Octavio Arizmendi, Jorge Garza-Vargas, and Daniel Perales.
\newblock Finite free cumulants: multiplicative convolutions, genus expansion
  and infinitesimal distributions.
\newblock {\em Transactions of the American Mathematical Society},
  376(06):4383--4420, 2023.

\bibitem{arizmendi2018cumulants}
Octavio Arizmendi and Daniel Perales.
\newblock Cumulants for finite free convolution.
\newblock {\em Journal of Combinatorial Theory, Series A}, 155:244--266, 2018.

\bibitem{bauerschmidt2020edge}
Roland Bauerschmidt, Jiaoyang Huang, Antti Knowles, and Horng-Tzer Yau.
\newblock Edge rigidity and universality of random regular graphs of
  intermediate degree.
\newblock {\em Geometric and functional analysis}, 30(3):693--769, 2020.

\bibitem{bilu2006lifts}
Yonatan Bilu and Nathan Linial.
\newblock Lifts, discrepancy and nearly optimal spectral gap.
\newblock {\em Combinatorica}, 26(5):495--519, 2006.

\bibitem{borcea2009lee}
Julius Borcea and Petter Br{\"a}nd{\'e}n.
\newblock The lee-yang and p{\'o}lya-schur programs. i. linear operators
  preserving stability.
\newblock {\em Inventiones mathematicae}, 177(3):541--569, 2009.

\bibitem{borcea2009lee2}
Julius Borcea and Petter Br{\"a}nd{\'e}n.
\newblock The lee-yang and p{\'o}lya-schur programs. ii. theory of stable
  polynomials and applications.
\newblock {\em Communications on Pure and Applied Mathematics: A Journal Issued
  by the Courant Institute of Mathematical Sciences}, 62(12):1595--1631, 2009.

\bibitem{bordenave2015new}
Charles Bordenave.
\newblock A new proof of friedman's second eigenvalue theorem and its extension
  to random lifts.
\newblock {\em arXiv preprint arXiv:1502.04482}, 2015.

\bibitem{bordenave2019eigenvalues}
Charles Bordenave and Beno{\^\i}t Collins.
\newblock Eigenvalues of random lifts and polynomials of random permutation
  matrices.
\newblock {\em Annals of Mathematics}, 190(3):811--875, 2019.

\bibitem{bownik2023akemann}
Marcin Bownik.
\newblock On akemann-weaver conjecture.
\newblock {\em arXiv preprint arXiv:2303.12954}, 2023.

\bibitem{chen2024new}
Chi-Fang Chen, Jorge Garza-Vargas, Joel~A Tropp, and Ramon van Handel.
\newblock A new approach to strong convergence.
\newblock {\em arXiv preprint arXiv:2405.16026}, 2024.

\bibitem{clark2006ramanujan}
Pete Clark.
\newblock Ramanujan graphs and shimura curves, 2006.

\bibitem{cohen2016ramanujan}
Michael~B Cohen.
\newblock Ramanujan graphs in polynomial time.
\newblock In {\em 2016 IEEE 57th Annual Symposium on Foundations of Computer
  Science (FOCS)}, pages 276--281. IEEE, 2016.

\bibitem{diaconis2004random}
Persi Diaconis and Alex Gamburd.
\newblock Random matrices, magic squares and matching polynomials.
\newblock {\em the electronic journal of combinatorics}, pages R2--R2, 2004.

\bibitem{diaconis1994eigenvalues}
Persi Diaconis and Mehrdad Shahshahani.
\newblock On the eigenvalues of random matrices.
\newblock {\em Journal of Applied Probability}, 31(A):49--62, 1994.

\bibitem{friedman1993some}
Joel Friedman.
\newblock Some geometric aspects of graphs and their eigenfunctions.
\newblock 1993.

\bibitem{friedman2003relative}
Joel Friedman.
\newblock Relative expanders or weakly relatively ramanujan graphs.
\newblock 2003.

\bibitem{friedman2008proof}
Joel Friedman.
\newblock {\em A proof of Alon's second eigenvalue conjecture and related
  problems}.
\newblock American Mathematical Soc., 2008.

\bibitem{friedman2019relativized}
Joel Friedman and David Kohler.
\newblock On the relativized alon second eigenvalue conjecture i: Main
  theorems, examples, and outline of proof.
\newblock {\em arXiv preprint arXiv:1911.05688}, 2019.

\bibitem{godsil1978matching}
Christopher~David Godsil and Ivan Gutman.
\newblock {\em On the matching polynomial of a graph}.
\newblock University of Melbourne Melbourne, 1978.

\bibitem{gorin2020crystallization}
Vadim Gorin and Adam~W Marcus.
\newblock Crystallization of random matrix orbits.
\newblock {\em International Mathematics Research Notices}, 2020(3):883--913,
  2020.

\bibitem{greenberg1995spectrum}
Yoseph Greenberg.
\newblock {\em On the spectrum of graphs and their universal covering}.
\newblock PhD thesis, Hebrew University, 1995.

\bibitem{gribinski2024theory}
Aurelien Gribinski.
\newblock A theory of singular values for finite free probability.
\newblock {\em Journal of Theoretical Probability}, 37(2):1257--1298, 2024.

\bibitem{hall2018ramanujan}
Chris Hall, Doron Puder, and William~F Sawin.
\newblock Ramanujan coverings of graphs.
\newblock {\em Advances in Mathematics}, 323:367--410, 2018.

\bibitem{heilmann1972theory}
Ole~J Heilmann and Elliott~H Lieb.
\newblock Theory of monomer-dimer systems.
\newblock {\em Communications in mathematical Physics}, 25(3):190--232, 1972.

\bibitem{hide2023near}
Will Hide and Michael Magee.
\newblock Near optimal spectral gaps for hyperbolic surfaces.
\newblock {\em Annals of Mathematics}, 198(2):791--824, 2023.

\bibitem{huang2024optimal}
Jiaoyang Huang, Theo McKenzie, and Horng-Tzer Yau.
\newblock Optimal eigenvalue rigidity of random regular graphs.
\newblock {\em arXiv preprint arXiv:2405.12161}, 2024.

\bibitem{huang2022edge}
Jiaoyang Huang and Horng-Tzer Yau.
\newblock Edge universality of sparse random matrices.
\newblock {\em arXiv preprint arXiv:2206.06580}, 2022.

\bibitem{keating2000random}
Jon~P Keating and Nina~C Snaith.
\newblock Random matrix theory and $\zeta$ (1/2+ it).
\newblock {\em Communications in Mathematical Physics}, 214:57--89, 2000.

\bibitem{lubotzky1994discrete}
Alex Lubotzky.
\newblock {\em Discrete groups, expanding graphs and invariant measures},
  volume 125.
\newblock Springer Science \& Business Media, 1994.

\bibitem{lubotzky1998not}
Alexander Lubotzky and Tatiana Nagnibeda.
\newblock Not every uniform tree covers ramanujan graphs.
\newblock {\em Journal of Combinatorial Theory, Series B}, 74(2):202--212,
  1998.

\bibitem{lubotzky1988ramanujan}
Alexander Lubotzky, Ralph Phillips, and Peter Sarnak.
\newblock Ramanujan graphs.
\newblock {\em Combinatorica}, 8(3):261--277, 1988.

\bibitem{marcus2016solution}
Adam Marcus and Nikhil Srivastava.
\newblock The solution of the kadison-singer problem.
\newblock {\em Current Developments in Mathematics}, 2016(1):111--143, 2016.

\bibitem{marcus2021polynomial}
Adam~W Marcus.
\newblock Polynomial convolutions and (finite) free probability.
\newblock {\em arXiv preprint arXiv:2108.07054}, 2021.

\bibitem{marcus2014ramanujan}
Adam~W Marcus, Daniel~A Spielman, and Nikhil Srivastava.
\newblock Ramanujan graphs and the solution of the kadison-singer problem.
\newblock In {\em 2014 International Congress of Mathematicans, ICM 2014},
  pages 363--386. KYUNG MOON SA Co. Ltd., 2014.

\bibitem{marcus2015interlacing}
Adam~W Marcus, Daniel~A Spielman, and Nikhil Srivastava.
\newblock Interlacing families i: Bipartite ramanujan graphs of all degrees.
\newblock {\em Annals of Mathematics}, pages 307--325, 2015.

\bibitem{if2}
Adam~W Marcus, Daniel~A Spielman, and Nikhil Srivastava.
\newblock Interlacing families ii: Mixed characteristic polynomials and the
  kadison—singer problem.
\newblock {\em Annals of Mathematics}, pages 327--350, 2015.

\bibitem{marcus2018interlacing}
Adam~W Marcus, Daniel~A Spielman, and Nikhil Srivastava.
\newblock Interlacing families iv: Bipartite ramanujan graphs of all sizes.
\newblock {\em SIAM Journal on computing}, 47(6):2488--2509, 2018.

\bibitem{marcus2022finite}
Adam~W Marcus, Daniel~A Spielman, and Nikhil Srivastava.
\newblock Finite free convolutions of polynomials.
\newblock {\em Probability Theory and Related Fields}, 182(3):807--848, 2022.

\bibitem{margulis1988explicit}
Grigorii~Aleksandrovich Margulis.
\newblock Explicit group-theoretical constructions of combinatorial schemes and
  their application to the design of expanders and concentrators.
\newblock {\em Problemy peredachi informatsii}, 24(1):51--60, 1988.

\bibitem{mckay1981expected}
Brendan~D McKay.
\newblock The expected eigenvalue distribution of a large regular graph.
\newblock {\em Linear Algebra and its applications}, 40:203--216, 1981.

\bibitem{miller2008distribution}
Steven~J Miller, Tim Novikoff, and Anthony Sabelli.
\newblock The distribution of the largest nontrivial eigenvalues in families of
  random regular graphs.
\newblock {\em Experimental Mathematics}, 17(2):231--244, 2008.

\bibitem{mohanty2019x}
Sidhanth Mohanty and Ryan O'Donnell.
\newblock $ x $-ramanujan graphs.
\newblock {\em arXiv preprint arXiv:1904.03500}, 2019.

\bibitem{puder2015expansion}
Doron Puder.
\newblock Expansion of random graphs: New proofs, new results.
\newblock {\em Inventiones mathematicae}, 201(3):845--908, 2015.

\bibitem{ravichandran2020mixed}
Mohan Ravichandran and Jonathan Leake.
\newblock Mixed determinants and the kadison--singer problem.
\newblock {\em Mathematische Annalen}, 377(1):511--541, 2020.

\bibitem{viennot2006heaps}
G{\'e}rard~Xavier Viennot.
\newblock Heaps of pieces, i: Basic definitions and combinatorial lemmas.
\newblock In {\em Combinatoire {\'e}num{\'e}rative: Proceedings of the
  “Colloque de combinatoire {\'e}num{\'e}rative”, held at Universit{\'e} du
  Qu{\'e}bec {\`a} Montr{\'e}al, May 28--June 1, 1985}, pages 321--350.
  Springer, 2006.

\bibitem{voiculescu1991limit}
Dan Voiculescu.
\newblock Limit laws for random matrices and free products.
\newblock {\em Inventiones mathematicae}, 104(1):201--220, 1991.

\bibitem{walsh1922location}
Joseph~L Walsh.
\newblock On the location of the roots of certain types of polynomials.
\newblock {\em Transactions of the American Mathematical Society},
  24(3):163--180, 1922.

\bibitem{wormald1999models}
Nicholas~C Wormald et~al.
\newblock Models of random regular graphs.
\newblock {\em London mathematical society lecture note series}, pages
  239--298, 1999.

\end{thebibliography}
\end{document}